# SEMIPARAMETRIC DENSITY ESTIMATION BY LOCAL $L_2$-FITTING[1]

By Kanta Naito

*Shimane University*

This article examines density estimation by combining a parametric approach with a nonparametric factor. The plug-in parametric estimator is seen as a crude estimator of the true density and is adjusted by a nonparametric factor. The nonparametric factor is derived by a criterion called local $L_2$-fitting. A class of estimators that have multiplicative adjustment is provided, including estimators proposed by several authors as special cases, and the asymptotic theories are developed. Theoretical comparison reveals that the estimators in this class are better than, or at least competitive with, the traditional kernel estimator in a broad class of densities. The asymptotically best estimator in this class can be obtained from the elegant feature of the bias function.

**1. Introduction.** Smoothing is a very important area of statistical analysis and has a wide range of applications in mathematical sciences. The present article is concerned especially with density estimation. Let $X_1, \dots, X_n$ be independently and identically distributed with density $f$. The problem is in estimating the density function $f$ from the data. In considering this problem, two approaches exist.

The first is called the parametric approach. In this approach, we prepare a parametric model

$$\{g(x,\theta) : \theta \in \Theta\},$$

where $\theta$ is a $p$-dimensional parameter vector and $\Theta$ is the parameter space in $\mathbb{R}^p$. In practice the family of densities is constructed from previous experience and preanalysis of the underlying structure. Then estimation of the

Received August 1998; revised April 2003.

[1]Supported in part by Grants-in-Aid for scientific research 11780171 and 13780174 from the Japan Society for the Promotion of Science.

density function is replaced by estimation of the unknown parameter vector $\theta$. Finally, we obtain a density estimator

$$\hat{f}(x) = g(x, \hat{\theta}),$$

where $\hat{\theta}$ is an estimator. This approach is called the plug-in parametric approach and is justified only when the true $f$ is exactly as in the model or at least in the neighborhood of the model.

The other approach is nonparametric. Several methods for nonparametric density estimation have been proposed and investigated. Izenman (1991) summarized a number of these methods. A representative method is the traditional kernel density estimator of $f$,

$$\tilde{f}(x) = \frac{1}{n}\sum_{i=1}^{n} K_h(X_i - x), \tag{1.1}$$

where $K_h(z) = h^{-1}K(h^{-1}z)$, $K(\cdot)$ is some chosen density which is symmetric about zero, and $h$ is the bandwidth. The basic properties of $\tilde{f}$ are well known and under smoothness conditions we have

$$\begin{aligned} E\tilde{f}(x) &= f(x) + \frac{h^2}{2}\mu_{2,K} f''(x) + O(h^4), \\ \operatorname{Var}\tilde{f}(x) &= \frac{R(K)}{nh} f(x) - \frac{f(x)^2}{n} + O\left(\frac{h}{n}\right), \end{aligned} \tag{1.2}$$

where $\mu_{\ell,G} = \int z^{\ell} G(z)\,dz$ and $R(G) = \int G(z)^2\,dz$ for some kernel function $G$ [cf. Simonoff (1996) and Wand and Jones (1995)]. The traditional kernel estimator is by construction completely nonparametric in the sense that it has no preferences and works reasonably well for almost all shapes of densities. Like the kernel estimator, all nonparametric methods can be used without the structural assumption that the underlying structure is controlled or captured by a finite-dimensional parameter. Thus, nonparametric approaches have attractive flexibility; however, the parametric model is difficult to discount because a well-estimated structure by the parametric approach is easy to understand.

This motivates us to propose an approach which includes both the parametric approach and the nonparametric approach. We propose and investigate a class of semiparametric density estimators which have precision comparable to, and sometime better than, that of $\tilde{f}$. One class considered herein is the set of density estimators derived from *the local $L_2$-fitting criterion with index $\alpha$*. In the proposed approach, the parametric plug-in density estimator $g(x, \hat{\theta})$ is utilized, but it is seen as a crude guess of $f(x)$. This initial parametric approximation is adjusted via multiplication by an adjustment factor $\xi = \xi(x)$. That is, the initial approximation is adjusted via the form $g(x,\hat{\theta})\xi$. The local fitting approach is used to determine the adjustment factor. Throughout the present article, $\xi = \xi(x)$ is determined by minimization of the empirical version of the function

$$Q(x,\xi|\alpha) = \int K_h(t-x)\frac{\{f(t) - g(x,\hat{\theta})\xi\}^2}{g(t,\hat{\theta})^{\alpha}}\,dt \tag{1.3}$$

for a fixed target point $x$. This method is called the local $L_2$-fitting criterion, where $\alpha$ is a real number called the index. Observe that local fitting is obtained using the kernel function $K$. The symmetric density $K$ creates the fitting locally around the target point $x$. This local approach is based on the simple intuition that observed data which are far from the target point $x$ do not have information about the adjustment. The minimizer of the empirical version of (1.3) is our objective and is denoted by $\hat{\xi} = \hat{\xi}(x)$. Using this $\hat{\xi}$, we finally obtain a density estimator $\hat{f}(x) = g(x,\hat{\theta})\hat{\xi}(x)$. This approach is shown to be effective and yields a theoretically good estimator in the sense of mean integrated squared error (MISE). A similar but somewhat different approach was proposed by Copas (1995) in conjunction with the likelihood method under censoring. Eguchi and Copas (1998) also discussed a class of local likelihood methods and developed asymptotics under a large bandwidth $h$. Their approach is the local estimation of $\theta$ in the model $g(x,\theta)$ and the adjustment factor $\xi$ does not appear. The present approach is the local estimation of $\xi$ using a previously obtained plug-in parametric estimator $g(x,\hat{\theta})$.

This multiplicative approach is closely related to studies performed by Hjort and Glad (1995) and Hjort and Jones (1996). Hjort and Glad (1995) proposed a density estimator based on the *naive* estimator of $\xi$. In addition, Hjort and Jones (1996) suggested and investigated two versions of multiplicative density estimators. One class of density estimators considered here includes these estimators as special cases, so this article may be seen as a generalization of these previous works.

The class of density estimators is developed in Section 2, and the estimators proposed by Hjort and Glad (1995) and Hjort and Jones (1996) are reviewed through examples. The behavior of the present estimators is investigated in Section 3, which also reveals that the present result is indeed a generalization of the results of Hjort and Glad (1995). The variance of the present estimator is the same as that of the traditional kernel estimator $\tilde{f}$, but the structure of the bias has a different form that depends on the initial parametric approximation. As an important property, we confirm that if $f$ is in the model, the estimator has reduced bias. Approximate or asymptotic MISE (AMISE) is derived in Section 4. Furthermore, the best estimator in the class is determined from the simple result that the bias is *linear* in $\alpha$. In Section 5 we compare the present estimator with $\tilde{f}$ for the case in which $f$ belongs to a class of normal mixture densities. In particular, a comparison is performed for 15 different test densities proposed by Marron and Wand (1992). In addition, a similar comparison for the case in which $f$ is the skew-normal distribution proposed by Azzalini (1985) is also discussed in Section 5. In Section 6 a simple algorithm to choose the best $\alpha$ is proposed. This algorithm is a variant of that used by Hjort and Glad (1995). Furthermore, two methods of data-based selection of $\alpha$ are discussed, and theoretical results and the practical algorithm are documented. These methods are constructed by reference to the theory of estimating the density functional discussed by Hall and Marron (1987) and Wand and Jones [(1995), Section 3.5]. Finite sample performance of the proposed estimators, and comparison to the $\tilde{f}$, the Hjort and Glad and the Hjort and Jones estimators are investigated by Monte Carlo simulation in Section 7. Supplementary remarks are presented in Section 8. It is trivial that the integral of the estimator is not unity, but the expansion formula as $h$ tends to zero shows that it is $1 + O(h^4)$ provided that we adopt a Gaussian density as an initial parametric model. A practical expression of the proposed estimator under the case using a Gaussian kernel and model is presented. Proofs of the theoretical results are presented in Section 9.

**2. Local $L_2$-fitting criterion.** This section is devoted to the construction of the present density estimator. First, we prepare a plug-in parametric density estimator $g(x, \hat{\theta})$, where $\hat{\theta}$ is an estimator of the least false value $\theta_0$ according to a certain distance measure between $f$ and $g(\cdot, \theta)$. The maximum likelihood estimator is a representative candidate for $\hat{\theta}$ in which the distance measure is known as the Kullback–Leibler distance $\int f(x) \log\{f(x)/g(x, \theta)\}\, dx$ and $\theta_0$ is defined as the minimizer of the Kullback–Leibler distance on $\theta$. This parametric estimator is seen as a crude guess of $f$. Next, we aim to adjust this initial approximation by the form $g(x, \hat{\theta})\xi$, where $\xi = \xi(x)$ is the adjustment factor. The problem is determination of $\xi$. To explain this method more clearly and to introduce the approaches proposed by Hjort and Glad (1995) and Hjort and Jones (1996), we present three examples below. Note that the kernel function $K$ is a symmetric density and the notation utilized in (1.1) and (1.2) is used throughout.

EXAMPLE 1 (Hjort and Jones estimator). To determine the adjustment factor $\xi$, Hjort and Jones (1996) suggested that the function of $\xi$ is

$$q(x, \xi) = \int K_h(t - x)\{f(t) - g(t, \hat{\theta})\xi\}^2\, dt.$$

The optimal $\xi$ is determined by minimization of the estimate of $q(x, \xi)$ on $\xi$. That is, we seek to minimize

$$q_n(x, \xi) = \xi^2 \int K_h(t - x) g(t, \hat{\theta})^2\, dt - \frac{2\xi}{n} \sum_{i=1}^{n} K_h(X_i - x) g(X_i, \hat{\theta}),$$

which gives

$$\hat{\xi} = \hat{\xi}(x) = \arg\min_{\xi} q_n(x,\xi) = \frac{n^{-1}\sum_{i=1}^{n} K_h(X_i - x) g(X_i, \hat{\theta})}{\int K_h(t-x) g(t,\hat{\theta})^2\, dt}.$$

The density estimator is obtained by

$$\hat{f}_{\mathrm{HJ}}(x) = g(x,\hat{\theta})\hat{\xi}(x) = g(x,\hat{\theta}) \frac{n^{-1}\sum_{i=1}^{n} K_h(X_i - x) g(X_i, \hat{\theta})}{\int K_h(t-x) g(t,\hat{\theta})^2\, dt}. \tag{2.1}$$

Although not fully discussed, this $\hat{f}_{\mathrm{HJ}}$ is the resultant estimator suggested by Hjort and Jones [(1996), page 1636].

Example 2 (Local likelihood estimator). The factor $\xi$ is determined by minimizing the empirical form of

$$\ell(x,\xi) = \int K_h(t-x)\left[f(t)\log\frac{f(t)}{g(t,\hat{\theta})\xi} - \{f(t) - g(t,\hat{\theta})\xi\}\right] dt,$$

which is equivalent to maximizing that of

$$L(x,\xi) = \int K_h(t-x)\{f(t)\log\{g(t,\hat{\theta})\xi\} - g(t,\hat{\theta})\xi\}\, dt.$$

The term $\ell(x,\xi)$ can be seen as a local version of the Kullback–Leibler distance from $f(x)$ to $g(x,\hat{\theta})\xi$. The resultant adjustment factor is

$$\hat{\xi} = \hat{\xi}(x) = \frac{\tilde{f}(x)}{\int K_h(t-x) g(t,\hat{\theta})\, dt}$$

and the ensuing estimator is

$$\begin{aligned} \hat{f}_{\mathrm{LL}}(x) &= g(x,\hat{\theta})\hat{\xi}(x) \\ &= g(x,\hat{\theta})\frac{\tilde{f}(x)}{\int K_h(t-x) g(t,\hat{\theta})\, dt} = \tilde{f}(x)\frac{g(x,\hat{\theta})}{\int K_h(t-x) g(t,\hat{\theta})\, dt}, \end{aligned} \tag{2.2}$$

where $\tilde{f}$ is as in (1.1). This $\hat{f}_{\mathrm{LL}}$ was proposed by Hjort and Jones [(1996), page 1635], who derived and discussed several estimators; $\hat{f}_{\mathrm{HJ}}$ and $\hat{f}_{\mathrm{LL}}$ are two special estimators with respect to the multiplicative adjustment scheme.

Example 3 (Hjort and Glad estimator). If we may assume $f(x) = g(x,\hat{\theta})\xi$, then true adjustment is $\xi = f(x)/g(x,\hat{\theta})$. Hjort and Glad (1995) proposed the *naive* estimator

$$\hat{\xi}(x) = \left(\widehat{f(x)/g(x,\hat{\theta})}\right) = \frac{1}{n}\sum_{i=1}^{n}\frac{K_h(X_i - x)}{g(X_i,\hat{\theta})},$$

which gives

$$\hat{f}_{\mathrm{HG}}(x) = g(x,\hat{\theta})\frac{1}{n}\sum_{i=1}^{n}\frac{K_h(X_i - x)}{g(X_i,\hat{\theta})}. \tag{2.3}$$

In Hjort and Glad (1995) the behavior of $\hat{f}_{\mathrm{HG}}$ was investigated and was shown to be better than the traditional kernel estimator in the sense of MISE on a certain class of normal mixture densities.

In the present article we are concerned with a function, namely (1.3), in conjunction with Examples 1–3. Considering the empirical version of $Q(x,\xi|\alpha)$ gives, by omitting the irrelevant term, the objective function

$$Q_n(x,\xi|\alpha) = \xi^2\int K_h(t-x)g(t,\hat{\theta})^{2-\alpha}\,dt - \frac{2\xi}{n}\sum_{i=1}^{n}K_h(X_i - x)g(X_i,\hat{\theta})^{1-\alpha}.$$

Obviously, $\alpha = 0$ gives $q_n(x,\xi)$, so $Q_n(x,\xi|\alpha)$ is a generalization of $q_n(x,\xi)$ in Example 1 and has weight function $g(t,\hat{\theta})^{-\alpha}$. The minimizer can be easily determined as

$$\hat{\xi} = \hat{\xi}(x) = \arg\min_{\xi} Q_n(x,\xi|\alpha) = \frac{n^{-1}\sum_{i=1}^{n}K_h(X_i - x)g(X_i,\hat{\theta})^{1-\alpha}}{\int K_h(t-x)g(t,\hat{\theta})^{2-\alpha}\,dt},$$

which is the proposed adjustment factor. Since the estimator depends on $\alpha$, by adding the symbol $\alpha$ we have

$$\hat{f}_{\alpha}(x) = g(x,\hat{\theta})\hat{\xi}(x) = g(x,\hat{\theta})\frac{n^{-1}\sum_{i=1}^{n}K_h(X_i - x)g(X_i,\hat{\theta})^{1-\alpha}}{\int K_h(t-x)g(t,\hat{\theta})^{2-\alpha}\,dt}. \tag{2.4}$$

From (2.1)–(2.4), the following relationships hold:

$$\hat{f}_0(x) = \hat{f}_{\mathrm{HJ}}(x), \qquad \hat{f}_1(x) = \hat{f}_{\mathrm{LL}}(x), \qquad \hat{f}_2(x) = \hat{f}_{\mathrm{HG}}(x).$$

The case $\alpha = 0$ is trivial. The case $\alpha = 1$ is confirmed by noting the definition $\ell(x,\xi)$ and the Taylor expansion of $(1+y)\log(1+y)$ at $y = 0$. This is also noted in Hjort and Jones (1996). The equality $\hat{f}_2 = \hat{f}_{\mathrm{HG}}$ claims that the naive estimator $\hat{\xi}$ proposed by Hjort and Glad (1995) is characterized by the minimizer of $Q_n(x,\xi|2)$. Therefore the estimators determined in Examples 1–3 are connected by $\alpha$. We thus propose a class of density estimators using $\alpha$ as the index. As described in the following sections, the introduction of $\alpha$ is essential and enables us to progress toward the theory of optimality in density estimation by the multiplicative adjustment scheme. In the following sections we discuss the behavior of estimators in this class. In addition, the best estimator in this class is determined.

**3. Asymptotic theory.** In this section, we investigate various statistically important quantities about $\hat{f}_\alpha$, such as bias and variance. From the features of $\hat{f}_\alpha$, it is trivial that its behavior depends on that of $\hat{\theta}$ included in the initial parametric approximation $g(x,\hat{\theta})$. To proceed with the theoretical study, we allow a somewhat more general setting for the choice of estimator $\hat{\theta}$. Let $F$ be the true distribution function, the cumulative of $f$, and let $F_n$ be the empirical distribution function. We consider functional estimators of $\theta$ of the form $\hat{\theta} = T(F_n)$ for some smooth functional $T$ having the influence function

$$I(x) = \lim_{\varepsilon \to 0} [T((1-\varepsilon)F + \varepsilon\delta_x) - T(F)]/\varepsilon,$$

where $\delta_x$ is the unit point mass at $x$, and assume that $\Sigma_I = E_f[I(X_i)I(X_i)^T]$ is finite. The best parametric approximation $g_0(x) = g(x,\theta_0)$ to $f(x)$ that $g(x,\hat{\theta})$ aims for is determined by $\theta_0 = T(F)$. It is well known for the case of the maximum likelihood estimator that $T(F)$ is defined as the solution of the equation $\int (\partial/\partial\theta) \log g(x,\theta)\, dF(x) = 0$, and so $I(x) = J^{-1}(\partial/\partial\theta)\log g(x,\theta_0)$, where $J = -E_f[(\partial^2/\partial\theta\,\partial\theta^T)\log g(X_i,\theta_0)]$. We may refer to Serfling (1980) for such a functional estimator. Under regularity conditions [see, e.g., Shao (1991)] we have

$$\hat{\theta} = \theta_0 + \frac{1}{n}\sum_{i=1}^{n} I(X_i) + \frac{d}{n} + \varepsilon_n, \tag{3.1}$$

where $\varepsilon_n = O_p(1/n)$ with mean $O(1/n^2)$. Then we have the following theorem.

THEOREM 1. *Let $g_0(x) = g(x,\theta_0)$, with $\theta_0 = T(F)$, be the best parametric approximation to $f$. Then, as $n \to \infty$, $h \to 0$,*

$$\operatorname{Bias} \hat{f}_\alpha(x) = \frac{h^2}{2}\mu_{2,K}\left[\frac{(g_0(x)^{1-\alpha}f(x))''}{g_0(x)^{1-\alpha}} - \frac{f(x)(g_0(x)^{2-\alpha})''}{g_0(x)^{2-\alpha}}\right] + O\left(h^4 + \frac{h^2}{n} + \frac{1}{n^2}\right),$$

$$\operatorname{Var} \hat{f}_\alpha(x) = \frac{R(K)}{nh} f(x) - \frac{f(x)^2}{n} + O\left(\frac{h}{n} + \frac{1}{n^2}\right).$$

The proof is included in Section 9. Note that the leading term of the variance of $\hat{f}_\alpha$ is independent of the estimation of $\theta$ and, with reference to (1.2), it is the same as that of $\tilde{f}$. Consistency of the density estimator requires both $h \to 0$ and $nh \to \infty$. The optimal size of $h$ is proportional to $n^{-1/5}$, which is also the same as that for $\tilde{f}$. Furthermore, it is worth noting that if $f$ is in the model $\{g(x,\theta):\theta\in\Theta\}$, that is, $g_0(x)=f(x)$, then the $O(h^2)$ term of the bias vanishes.

From the above observations, the essential difference between the behavior of $\hat{f}_\alpha$ and that of $\tilde{f}$ appears in the bias. As seen in the next section, the $O(h^2)$ term of the bias of $\hat{f}_\alpha$ has a nice expression (4.5), which allows the best estimator in the sense of MISE to be determined.

**4. Goodness of estimators.** In this section, the goodness of estimators is evaluated in the sense of MISE. In addition, $\hat{f}_\alpha$ and $\tilde{f}$ are compared. Let $\mathcal{R}(\bar{f})$ denote the integral of the squared $O(h^2)$ term of the bias of a density estimator $\bar{f}$. From Theorem 1 and (1.2), the AMISE of $\hat{f}_\alpha$ and $\tilde{f}$ are, respectively, given by

$$\mathrm{AMISE}(\hat{f}_\alpha)=\frac{h^4}{4}\mu_{2,K}^2\mathcal{R}(\hat{f}_\alpha)+\frac{R(K)}{nh}$$

and

$$\mathrm{AMISE}(\tilde{f})=\frac{h^4}{4}\mu_{2,K}^2\mathcal{R}(\tilde{f})+\frac{R(K)}{nh},$$

where

$$\mathcal{R}(\hat{f}_\alpha)=\int\left[\frac{(g_0(x)^{1-\alpha}f(x))''}{g_0(x)^{1-\alpha}}-\frac{f(x)(g_0(x)^{2-\alpha})''}{g_0(x)^{2-\alpha}}\right]^2dx, \tag{4.1}$$

$$\mathcal{R}(\tilde{f})=\int\{f''(x)\}^2\,dx. \tag{4.2}$$

So it suffices to compare $\mathcal{R}(\hat{f}_\alpha)$ and $\mathcal{R}(\tilde{f})$ in the AMISE comparison, provided that we use the same kernel function $K$. The AMISE comparison will be discussed for special choices of the underlying $f$, using the same kernel.

Now we consider the function in the bracket in (4.1) to discover the best estimator. Let us define

$$b_1(x)=f''(x)-f(x)\frac{g_0{}''(x)}{g_0(x)}, \tag{4.3}$$

$$b_2(x)=2\left\{\frac{g_0{}'(x)f'(x)}{g_0(x)}-f(x)\left(\frac{g_0{}'(x)}{g_0(x)}\right)^2\right\}. \tag{4.4}$$

Then it is easily verified that

$$\frac{(g_0(x)^{1-\alpha}f(x))''}{g_0(x)^{1-\alpha}}-\frac{f(x)(g_0(x)^{2-\alpha})''}{g_0(x)^{2-\alpha}}=\{b_1(x)+b_2(x)\}-\alpha b_2(x). \tag{4.5}$$

That is, the $O(h^2)$ term of the bias of $\hat{f}_\alpha$ is *linear* in $\alpha$. Therefore, writing

$$c_1=\int\{b_2(x)\}^2\,dx, \tag{4.6}$$

$$c_2 = \int b_2(x)\{b_1(x) + b_2(x)\}\, dx, \tag{4.7}$$

$$c_3 = \int \{b_1(x) + b_2(x)\}^2\, dx, \tag{4.8}$$

we obtain

$$\mathcal{R}(\hat{f}_\alpha) = c_1\alpha^2 - 2c_2\alpha + c_3. \tag{4.9}$$

Using (4.9), we have the leading terms of the integrated squared bias of $\hat{f}_{\mathrm{HJ}}$, $\hat{f}_{\mathrm{LL}}$ and $\hat{f}_{\mathrm{HG}}$ by substituting $\alpha = 0, 1$ and 2, respectively. For instance, $c_3 = \mathcal{R}(\hat{f}_0)$ is found to be the integrated squared bias of $\hat{f}_{\mathrm{HJ}}$. The quadratic expression of (4.9) establishes the following proposition.

PROPOSITION 1. *$\mathcal{R}(\hat{f}_\alpha)$ is minimized over $\alpha$ at*

$$\alpha_o = \frac{c_2}{c_1} \tag{4.10}$$

*and its minimum value is*

$$\min \mathcal{R}(\hat{f}_\alpha) = c_3 - \frac{(c_2)^2}{c_1}, \tag{4.11}$$

*where $c_1$–$c_3$ are given in* (4.6)–(4.8), *respectively.*

The linear structure (4.5) is essential in the derivation of Proposition 1. This is obtained by introducing $\alpha$ through the weighting $g(t, \hat{\theta})^{-\alpha}$ in $Q(x, \xi|\alpha)$, so that such a generalization indeed has an advantage. Theoretically, the *ideal* estimator $\hat{f}_{\alpha_o}$ is the best estimator in the class which surpasses estimators $\hat{f}_{\mathrm{HJ}}$, $\hat{f}_{\mathrm{LL}}$ and $\hat{f}_{\mathrm{HG}}$ in the sense of AMISE.

## 5. Asymptotic comparison.

In this section the proposed $\hat{f}_\alpha$ is compared to $\tilde{f}$ based on the AMISE formulas described in Section 4.

### 5.1. *Comparison in normal mixture.*

Here we compare $\hat{f}_\alpha$ and $\tilde{f}$ for the case in which $f$ belongs to the class of normal mixture densities. Let

$$f(x) = \sum_{i=1}^{k} p_i f_i(x),$$

where

$$f_i(x) = \frac{1}{\sigma_i}\phi\left(\frac{x - \mu_i}{\sigma_i}\right) \equiv \phi_{\sigma_i}(x - \mu_i),$$

$\phi$ is the standard normal density function and $\sum_{i=1}^{k} p_i = 1$. The family of such mixtures forms a very wide and flexible class of densities. Marron and Wand (1992) studied such mixtures and singled out 15 different densities which are often used as test densities in the study of the performance of density estimators [Hjort and Glad (1995), Jones and Signorini (1997) and Jones, Linton and Nielsen (1995)]. It is easy to see that

$$\mu_0 \equiv \int x f(x)\, dx = \sum_{i=1}^{k} p_i \mu_i,$$

$$\sigma_0^2 \equiv \int (x-\mu_0)^2 f(x)\, dx = \sum_{i=1}^{k} p_i \{\sigma_i^2 + (\mu_i - \mu_0)^2\}.$$

For the present estimator $\hat{f}_\alpha$, we adopt here the normal density $\phi_{\sigma_0}(x - \mu_0)$ as $g_0(x) = g(x, \theta_0)$. This corresponds to the use of maximum likelihood estimates (MLE) for estimation of $\theta_0$, since the normal density that has mean $\mu_0$ and variance $\sigma_0^2$ minimizes the Kullback–Leibler distance from $f(x)$ to $g(x,\theta) = \phi_\sigma(x-\mu)$, where $\theta = (\mu, \sigma^2)$ and $\theta_0 = (\mu_0, \sigma_0^2)$.

The previous section indicates that the AMISE comparison is performed by comparing $\mathcal{R}(\hat{f}_\alpha)$ and $\mathcal{R}(\tilde{f})$. Both can be calculated through (4.1) and (4.2) using numerical integration. However, when $f$ is a normal mixture and $g_0$ is normal, we obtain the analytic expression of $\mathcal{R}(\hat{f}_\alpha)$ by obtaining those of $c_1$, $c_2$ and $c_3$. Referring to (4.3) and (4.4), direct computation yields

$$b_1(x) = \sum_{i=1}^{k} p_i f_i(x) \left\{ \frac{1}{\sigma_i^2} H_2\left(\frac{x-\mu_i}{\sigma_i}\right) - \frac{1}{\sigma_0^2} H_2\left(\frac{x-\mu_0}{\sigma_0}\right) \right\},$$

$$b_2(x) = 2\sum_{i=1}^{k} p_i f_i(x) \left\{ \frac{1}{\sigma_0\sigma_i} H_1\left(\frac{x-\mu_0}{\sigma_0}\right) H_1\left(\frac{x-\mu_i}{\sigma_i}\right) - \frac{1}{\sigma_0^2} H_1^2\left(\frac{x-\mu_0}{\sigma_0}\right) \right\},$$

where $H_k$ is the $k$th order Hermite polynomial. Since $c_1$, $c_2$ and $c_3$ are all integrals of these functions, we find their analytic expressions using the properties of the Hermite polynomials. The detailed calculations are found in Naito [(1998), Sections 4 and 6]. On the other hand, the expression of $\mathcal{R}(\tilde{f})$ has already been presented in Marron and Wand (1992). Thus, by using (4.1) and (4.2), we compare $\hat{f}_\alpha$ and $\tilde{f}$ for 15 representative test densities used in Marron and Wand (1992). The values of the ratio $\mathcal{R}(\hat{f}_\alpha)/\mathcal{R}(\tilde{f})$ for $\alpha = 0, 1, 2, \alpha_o$ are tabulated in Table 1, in which the case number corresponds to that used in Marron and Wand (1992). The entries in column $\alpha_o$ are the values of $\alpha_o$ for each case. Since #1 is normal, $\mathcal{R}(\hat{f}_\alpha) = 0$ for all $\alpha$, so that the ratio is always zero in the #1 row. For example, in #6, which corresponds to a bimodal density, the value of $\mathcal{R}(\hat{f}_0)/\mathcal{R}(\tilde{f})$ is 1.7434 and that of $\mathcal{R}(\hat{f}_2)/\mathcal{R}(\tilde{f})$ is 0.7705, and for #6, the minimum of the ratio is attained at $\alpha_o = 1.9394$ and its minimum value is 0.7696.

TABLE 1
*Comparison in normal mixture*[a]

| f | $\alpha = 0$ | $\alpha = 1$ | $\alpha = 2$ | $\alpha = \alpha_o$ | $\alpha_o$ |
|---|---|---|---|---|---|
| #1 | 0.0000 | 0.0000 | 0.0000 | 0.0000 | — |
| #2 | 1.0448 | 0.3947 | 0.2460 | 0.2356 | 1.7968 |
| #3 | 1.0239 | 0.9986 | 0.9925 | 0.9922 | 1.8207 |
| #4 | 1.0010 | 0.9799 | 0.9606 | 0.8719 | 11.7075 |
| #5 | 1.0436 | 0.8826 | 0.7822 | 0.7414 | 3.1606 |
| #6 | 1.7434 | 0.9980 | 0.7705 | 0.7696 | 1.9394 |
| #7 | 1.4821 | 0.9829 | 0.8524 | 0.8485 | 1.8541 |
| #8 | 1.5398 | 1.0114 | 0.9007 | 0.8892 | 1.7651 |
| #9 | 1.3088 | 1.0010 | 0.9178 | 0.9159 | 1.8706 |
| #10 | 1.0512 | 0.9947 | 0.9791 | 0.9788 | 1.8787 |
| #11 | 1.0003 | 1.0000 | 0.9999 | 0.9999 | 1.8597 |
| #12 | 1.0236 | 1.0036 | 1.0025 | 1.0007 | 1.5589 |
| #13 | 1.0005 | 1.0000 | 0.9999 | 0.9999 | 1.7840 |
| #14 | 1.0030 | 1.0004 | 1.0002 | 1.0000 | 1.5897 |
| #15 | 1.0127 | 1.0013 | 1.0001 | 0.9994 | 1.6190 |

[a]Values of the ratio $\mathcal{R}(\hat{f}_\alpha)/\mathcal{R}(\tilde{f})$ are tabulated for the 15 densities in Marron and Wand (1992). Values of the optimal index $\alpha_o$ defined in (4.10) are listed in the $\alpha_o$ column for each case.

We can confirm that Proposition 1 holds and $\hat{f}_{\alpha_o}$ is better than, or at least competitive with, $\tilde{f}$ for all cases in this comparison. Furthermore, it is worth noting that $\alpha_o$ is around 2, except for #4 and #5. This reveals that the Hjort and Glad estimator $\hat{f}_{\mathrm{HG}} = \hat{f}_2$ is also good for almost all cases.

5.2. *Comparison in skew-normal.* Similar to the previous section, the comparison of $\hat{f}_\alpha$ and $\tilde{f}$ is performed for the case in which $f$ belongs to a class of skew-normal distributions discussed in Azzalini (1985). If a random variable $X$ has density $f(x) = 2\phi(x)\Phi(\lambda x)$, where $\Phi$ is the distribution function of the standard normal, then we say that $X$ has skew-normal distribution with parameter $\lambda$ and we denote this by $X \sim \mathrm{SN}(\lambda)$. Here SN(0) corresponds to the standard normal. We obtain from direct calculations that

$$f'(x) = 2\phi(x)s_1(x,\lambda), \qquad f''(x) = 2\phi(x)s_2(x,\lambda), \tag{5.1}$$

where

$$s_1(x,\lambda) = \lambda\phi(\lambda x) - H_1(x)\Phi(\lambda x),$$
$$s_2(x,\lambda) = H_2(x)\Phi(\lambda x) - (\lambda^3 + 2\lambda)H_1(x)\phi(\lambda x)$$

and $H_k$ is the $k$th order Hermite polynomial. In addition, we adopt the normal density as an initial approximation and the MLE for estimation of the parameter included in the parametric model. We have for $X \sim \mathrm{SN}(\lambda)$,

$$\mu_0 \equiv \int x f(x)\,dx = \sqrt{\frac{2}{\pi}}\frac{\lambda}{\sqrt{1+\lambda^2}},$$

$$\sigma_0^2 \equiv \int (x-\mu_0)^2 f(x)\,dx = 1 - \frac{2\lambda^2}{\pi(1+\lambda^2)},$$

which gives the least false parameter vector $\theta_0 = (\mu_0, \sigma_0^2)$ for $g_0(x) = \phi_{\sigma_0}(x-\mu_0)$. To find the best estimator, it is required to obtain $b_1(x)$ and $b_2(x)$ in (4.3) and (4.4), respectively. Direct computations yield

$$b_1(x) = 2\phi(x)\left[s_2(x,\lambda) - \frac{1}{\sigma_0^2}H_2\left(\frac{x-\mu_0}{\sigma_0}\right)\Phi(\lambda x)\right],$$

$$b_2(x) = -4\phi(x)\left[\frac{1}{\sigma_0}s_1(x,\lambda)H_1\left(\frac{x-\mu_0}{\sigma_0}\right) + \frac{1}{\sigma_0^2}H_1\left(\frac{x-\mu_0}{\sigma_0}\right)^2\Phi(\lambda x)\right].$$

Using these, we can obtain $\mathcal{R}(\hat{f}_\alpha)$, and we have from (4.2) and (5.1) that

$$\mathcal{R}(\tilde{f}) = \int \{2\phi(x)s_2(x,\lambda)\}^2\,dx.$$

Table 2 exhibits the comparison for $\lambda = 0(1)5$. For each $\lambda$ the ratio $\mathcal{R}(\hat{f}_\alpha)/\mathcal{R}(\tilde{f})$ is tabulated. Since $\lambda = 0$ implies $f = g_0$, the ratios are zero for all $\alpha$. For any $\lambda$ utilized in this comparison, we observe $\hat{f}_\alpha$ for $\alpha = 1, 2, \alpha_o$ are all superior to $\tilde{f}$.

TABLE 2
*Comparison in skew-normal*[a]

| f | $\alpha = 0$ | $\alpha = 1$ | $\alpha = 2$ | $\alpha = \alpha_o$ | $\alpha_o$ |
|---|---|---|---|---|---|
| $\lambda = 0$ | 0.0000 | 0.0000 | 0.0000 | 0.0000 | — |
| $\lambda = 1$ | 0.0762 | 0.0232 | 0.0134 | 0.0118 | 1.7270 |
| $\lambda = 2$ | 0.7636 | 0.2669 | 0.1645 | 0.1531 | 1.7594 |
| $\lambda = 3$ | 1.4625 | 0.5783 | 0.3945 | 0.3748 | 1.7624 |
| $\lambda = 4$ | 1.7888 | 0.7836 | 0.5839 | 0.5583 | 1.7480 |
| $\lambda = 5$ | 1.8678 | 0.8963 | 0.7133 | 0.6850 | 1.7320 |

[a]Values of the ratio $\mathcal{R}(\hat{f}_\alpha)/\mathcal{R}(\tilde{f})$ are tabulated for $\lambda = 0(1)5$ in $\mathrm{SN}(\lambda)$ proposed by Azzalini (1985). Values of the optimal index $\alpha_o$ defined in (4.10) are listed in the $\alpha_o$ column for each $\mathrm{SN}(\lambda)$.

**6. Index selection.** In this section, three data-based methods used to select the index $\alpha$ are discussed. These methods are somewhat intuitive, but the density estimators with the index obtained through these methods perform well, as shown in the simulation report in Section 7.

6.1. *Direct method.* We propose a data-based selection of $\alpha$ which is a derivative of that of $h$ discussed in Hjort and Glad [(1995), Section 6]. We consider the Hermite expansion given as

$$f(x) = \phi\left(\frac{x-\mu}{\sigma}\right)\frac{1}{\sigma}\left\{1 + \sum_{k=3}^{m}\frac{\gamma_k}{k!}H_k\left(\frac{x-\mu}{\sigma}\right)\right\}, \tag{6.1}$$

where $\gamma_0 = 1$ and $\gamma_1 = \gamma_2 = 0$. We know that $\gamma_k = E[H_k((X-\mu)/\sigma)]$. Simple but somewhat tedious computations, along with the Gaussian initial approximation $g_0(x) = \phi_\sigma(x-\mu)$ and $m = 5$, yield

$$c_1 = \frac{1}{\sigma^5\sqrt{\pi}}\left[\gamma_3^2\left(\frac{7}{16}\right) + \frac{\gamma_4^2}{9}\left(\frac{33}{32}\right) + \frac{\gamma_5^2}{144}\left(\frac{225}{64}\right) - \frac{\gamma_3\gamma_5}{6}\left(\frac{21}{32}\right)\right], \tag{6.2}$$

$$c_2 = \frac{1}{\sigma^5\sqrt{\pi}}\left[\gamma_3^2\left(\frac{3}{4}\right) + \frac{\gamma_4^2}{9}\left(\frac{32}{57}\right) + \frac{\gamma_5^2}{144}\left(\frac{195}{32}\right) - \frac{\gamma_3\gamma_5}{6}\left(\frac{39}{32}\right)\right], \tag{6.3}$$

$$c_3 = \frac{1}{\sigma^5\sqrt{\pi}}\left[\gamma_3^2\left(\frac{3}{2}\right) + \frac{\gamma_4^2}{9}\left(\frac{123}{32}\right) + \frac{\gamma_5^2}{144}\left(\frac{225}{16}\right) - \frac{\gamma_3\gamma_5}{2}\right]. \tag{6.4}$$

Here $c_i$, $i = 1, 2, 3$, are estimated in the usual manner by substituting

$$\hat{\gamma}_k = \frac{1}{n}\sum_{i=1}^{n} H_k\left(\frac{X_i - \hat{\mu}}{\hat{\sigma}}\right)$$

for $\gamma_k$, where $k = 3, 4, 5$, and by substituting $\hat{\sigma}$ for $\sigma$. The next step is to use nonparametric estimators of $c_1$ and $c_2$ defined by

$$\hat{c}_1(h) = \int \{\hat{b}_2(x;h)\}^2\,dx,$$

$$\hat{c}_2(h) = \int \hat{b}_2(x;h)\{\hat{b}_1(x;h) + \hat{b}_2(x;h)\}\,dx,$$

where

$$\hat{b}_1(x;h) = \frac{1}{n}\sum_{i=1}^{n}\left[\frac{1}{h^3}K''\left(\frac{x-X_i}{h}\right) - \frac{1}{h}K\left(\frac{x-X_i}{h}\right)\frac{g''(x,\hat{\theta})}{g(x,\hat{\theta})}\right],$$

$$\hat{b}_2(x;h) = \frac{2}{n}\sum_{i=1}^{n}\left[\frac{1}{h^2}K'\left(\frac{x-X_i}{h}\right)\frac{g'(x,\hat{\theta})}{g(x,\hat{\theta})} - \frac{1}{h}K\left(\frac{x-X_i}{h}\right)\left(\frac{g'(x,\hat{\theta})}{g(x,\hat{\theta})}\right)^2\right],$$

$K$ is a kernel, which may be different from that used in $\hat{f}_\alpha$, and $h$ is the bandwidth. Using these quantities, we choose $\alpha$ as follows. First, we obtain $\bar{c}_i$, $i = 1, 2, 3$, from (6.2)–(6.4), respectively, using $\hat{\gamma}_k$, $k = 3, 4, 5$, and $\hat{\sigma}$, under the assumption that the underlying distribution is approximated by the Hermite expansion. Then, referring to (4.11), $\mathcal{R}(\hat{f}_{\alpha_o})$ is estimated as

$$\bar{\mathcal{R}}(\hat{f}_{\alpha_o}) = \bar{c}_3 - \frac{(\bar{c}_2)^2}{\bar{c}_1}.$$

This gives a bandwidth

$$\bar{h} = \left\{ \frac{R(K)}{\mu_{2,K}^2} \right\}^{1/5} \bar{\mathcal{R}}(\hat{f}_{\alpha_o})^{-1/5} n^{-1/5}, \tag{6.5}$$

from which we have an estimate of the optimal index,

$$\hat{\alpha}_o^{[1]} = \frac{\hat{c}_2(\bar{h})}{\hat{c}_1(\bar{h})}. \tag{6.6}$$

6.2. *Two methods based on functional estimation.* Here we propose two methods based on estimation of the functional of $f$ and $g(x, \hat{\theta})$. Define

$$q_1(x) = \frac{g_0'(x)}{g_0(x)}, \tag{6.7}$$

$$q_2(x) = \frac{g_0''(x)}{g_0(x)} = q_1'(x) + \{q_1(x)\}^2, \tag{6.8}$$

where $g_0(x) = g(x, \theta_0)$. Using this notation, we have

$$c_1 = 4\int f'(x)^2 q_1(x)^2\,dx + 4\int f(x)^2 q_1(x)^4\,dx - 8\int f(x) f'(x) q_1(x)^3\,dx,$$

$$\begin{aligned} c_2 = c_1 &+ 2\int f'(x) f''(x) q_1(x)\,dx - 2\int f(x) f'(x) q_1(x) q_2(x)\,dx \\ &- 2\int f(x) f''(x) q_1(x)^2\,dx + 2\int f(x)^2 q_1(x)^2 q_2(x)\,dx. \end{aligned}$$

Under the sufficient smoothness condition for $f$, it follows that

$$\begin{aligned} &\int f'(x)^2 q_1(x)^2\,dx \\ &\quad = -E_f[f''(X) q_1(X)^2] - 2E_f[f'(X) q_1(X) q_2(X)] + 2E_f[f'(X) q_1(X)^3] \end{aligned}$$

and

$$\begin{aligned} &\int f'(x) f''(x) q_1(x)\,dx \\ &\quad = -E_f[f'''(X) q_1(X)] - E_f[f''(X) q_2(X)] + E_f[f''(X) q_1(X)^2]. \end{aligned}$$

These calculations allow us to define

$$\psi(p|r, s) \equiv E_f[f^{(p)}(X) q_1(X)^r q_2(X)^s]$$

for integers $p = 0, 1, 2, 3$, $r = 0, 1, 2, 3, 4$ and $s = 0, 1, 2$, where $f^{(p)}(x) = (d^p/dx^p) f(x)$ and $f^{(0)}(x) = f(x)$. Then we have

$$c_1 = 4\{\psi(0|4, 0) - \psi(2|2, 0) - 2\psi(1|1, 1)\},$$

$$c_2 = c_1 + 2\{\psi(0|2, 1) - \psi(3|1, 0) - \psi(2|0, 1) - \psi(1|1, 1)\},$$

so that the optimal $\alpha_o$ in (4.10) can be written in terms of $\psi$ as

$$\alpha_o = \frac{c_2}{c_1} = 1 + \frac{1}{2}\left[\frac{\psi(0|2,1) - \psi(3|1,0) - \psi(2|0,1) - \psi(1|1,1)}{\psi(0|4,0) - \psi(2|2,0) - 2\psi(1|1,1)}\right]$$
$$\equiv 1 + \frac{1}{2}\frac{\mathcal{N}}{\mathcal{D}},$$

where

$$\mathcal{N} = \psi(0|2,1) - \psi(3|1,0) - \psi(2|0,1) - \psi(1|1,1),$$
$$\mathcal{D} = \psi(0|4,0) - \psi(2|2,0) - 2\psi(1|1,1).$$

By the above reductions, data-based selection of $\alpha$ is accomplished by using an estimator of $\alpha_o$ defined by

$$\hat{\alpha}_o(g) = 1 + \frac{1}{2}\frac{\widehat{\mathcal{N}}_g}{\widehat{\mathcal{D}}_g}$$
$$= 1 + \frac{1}{2}\left[\frac{\hat{\psi}_g(0|2,1) - \hat{\psi}_g(3|1,0) - \hat{\psi}_g(2|0,1) - \hat{\psi}_g(1|1,1)}{\hat{\psi}_g(0|4,0) - \hat{\psi}_g(2|2,0) - 2\hat{\psi}_g(1|1,1)}\right],$$

where

$$\hat{\psi}_g(p|r,s) = \frac{1}{n(n-1)}\sum_{i\neq j}\hat{q}_1(X_i)^r\hat{q}_2(X_i)^s L_g^{(p)}(X_i - X_j)$$

is a nonparametric estimator of $\psi(p|r,s)$ that has a symmetric kernel $L$ and bandwidth $g$ that are possibly different from $K$ and $h$, respectively. In addition, $\hat{q}_1$ and $\hat{q}_2$ are, respectively, those of (6.7) and (6.8) using $g(x,\hat{\theta})$ rather than $g_0(x)$.

The behavior of $\hat{\alpha}_o(g)$ can be investigated by a method based on the theory of estimating the density functional [see, e.g., Section 3.5 in Wand and Jones (1995)]. Mean squared error (MSE) is adopted to evaluate $\widehat{\mathcal{N}}_g$ and $\widehat{\mathcal{D}}_g$, while $\hat{\alpha}_o(g)$ is evaluated by mean squared relative error (MSRE). Somewhat tedious calculations yield the following theorem:

THEOREM 2. *As $n \to \infty$ and $g \to 0$,*

$$\text{MSE}[\widehat{\mathcal{N}}_g] = \frac{g^4}{4}\mu_{2,L}^2\mathcal{N}[2]^2 + \frac{1}{2n^2g^5}\iint \lambda_{2|3}(x,z)^2\,dx\,dz + O(n^{-1}) + o(g^4 + n^{-2}g^{-5}), \tag{6.9}$$

$$\text{MSE}[\widehat{\mathcal{D}}_g] = \frac{g^4}{4}\mu_{2,L}^2\mathcal{D}[2]^2 + \frac{1}{2n^2g^5}\iint \kappa_2(x,z)^2\,dx\,dz + O(n^{-1}) + o(g^4 + n^{-2}g^{-5}), \tag{6.10}$$

$$\begin{aligned}\mathrm{MSRE}[\hat{\alpha}_o(g)] &= \frac{g^4}{16}\mu_{2,L}^2\left[\frac{\mathcal{N}[2]}{\mathcal{N}[0]} - \frac{\mathcal{D}[2]}{\mathcal{D}[0]}\right]^2 \\ &\quad + \frac{1}{8n^2g^5}\iint\left[\frac{\lambda_{2|3}(x,z)}{\mathcal{N}[0]} - \frac{\kappa_2(x,z)}{\mathcal{D}[0]}\right]^2 dx\,dz \\ &\quad + O(n^{-1}) + o(g^4 + n^{-2}g^{-5} + n^{-1}),\end{aligned} \tag{6.11}$$

*where*

$$\lambda_{p_2|p_1}(x,z) = f(x)[\{2L^{(p_2)}(z) + zL^{(p_1)}(z)\}q_2(x) - zL^{(p_1)}(z)q_1(x)^2],$$

$$\kappa_{p_2}(x,z) = f(x)[2L^{(p_2)}(z)q_1(x)^2]$$

*for even* $p_2$ *and* $p_1 = p_2 + 1$, *and*

$$\mathcal{N}[p] = \psi(p|2,1) - \psi(p+3|1,0) - \psi(p+2|0,1) - \psi(p+1|1,1)$$

$$\mathcal{D}[p] = \psi(p|4,0) - \psi(p+2|2,0) - 2\psi(p+1|1,1)$$

*for even* $p$ *with* $\mathcal{N}[0] = \mathcal{N}$ *and* $\mathcal{D}[0] = \mathcal{D}$.

The proof of Theorem 2 is presented in Section 9. From Theorem 2 the approximate mean squared error (AMSE)-optimal bandwidths for $\widehat{\mathcal{N}}_g$ and $\widehat{\mathcal{D}}_g$, and the approximate mean squared relative error (AMSRE)-optimal bandwidth for $\hat{\alpha}_o(g)$ are, respectively, given as

$$g_{\mathcal{N}}\text{-AMSE} = \left[\left(\frac{5}{2}\right)\frac{\iint \lambda_{2|3}(x,z)^2\,dx\,dz}{\mu_{2,L}^2\mathcal{N}[2]}\right]^{1/9} n^{-2/9},$$

$$g_{\mathcal{D}}\text{-AMSE} = \left[\left(\frac{5}{2}\right)\frac{\iint \kappa_2(x,z)^2\,dx\,dz}{\mu_{2,L}^2\mathcal{D}[2]}\right]^{1/9} n^{-2/9}$$

and

$$g_{\mathrm{AMSRE}} = \left[\left(\frac{5}{2}\right)\frac{\iint\{\mathcal{D}\lambda_{2|3}(x,z) - \mathcal{N}\kappa_2(x,z)\}^2\,dx\,dz}{\mu_{2,L}^2\{\mathcal{D}\mathcal{N}[2] - \mathcal{N}\mathcal{D}[2]\}^2}\right]^{1/9} n^{-2/9}.$$

Unfortunately, these bandwidths have the same defect as the plug-in method for bandwidth selection of the kernel density estimator: all of these bandwidths depend on unknown $\mathcal{N}[2]$, $\mathcal{D}[2]$, $\mathcal{N}$ and $\mathcal{D}$. Estimation of $\mathcal{N}[2]$ and $\mathcal{D}[2]$ is possible; however, their optimal bandwidths depend on $\mathcal{N}[4]$ and $\mathcal{D}[4]$. Furthermore, it can easily be recognized that this problem does not go away.

To overcome this problem, we utilize a simple estimate based on the Hermite expansion of (6.1). Equation (6.1) yields a pilot estimate of $f^{(p)}(x)$ as

$$\tilde{f}^{(p)}(x) = \frac{(-1)^p}{\hat{\sigma}^{p+1}}\phi\left(\frac{x-\hat{\mu}}{\hat{\sigma}}\right)\sum_{k=1}^{m}\frac{\hat{\gamma}_k}{k!}H_{k+p}\left(\frac{x-\hat{\mu}}{\hat{\sigma}}\right),$$

from which we have $\widetilde{\mathcal{N}}[6] = \tilde{\psi}(6|2,1) - \tilde{\psi}(9|1,0) - \tilde{\psi}(8|0,1) - \tilde{\psi}(7|1,1)$ as an estimate of $\mathcal{N}[6]$ using the component defined by

$$\tilde{\psi}(p|r,s) = \frac{(-1)^p}{\hat{\sigma}^p} \sum_{k=0}^{m} \frac{\hat{\gamma}_k}{k!} \left[ \frac{1}{n} \sum_{i=1}^{n} \frac{1}{\hat{\sigma}} \phi\left(\frac{X_i - \hat{\mu}}{\hat{\sigma}}\right) H_{k+p}\left(\frac{X_i - \hat{\mu}}{\hat{\sigma}}\right) \hat{q}_1(X_i)^r \hat{q}_2(X_i)^s \right].$$

An estimate $\widetilde{\mathcal{D}}[6]$ of $\mathcal{D}[6]$ can be obtained in the same manner.

In the following text we describe the algorithm used to obtain two estimates of $\alpha_o$. The notation utilized is

$$L^{[1]}(p_1,p_2) = \mu_{2,L^{(p_1)}L^{(p_1)}} + 4\mu_{0,L^{(p_2)}L^{(p_2)}} + 4\mu_{1,L^{(p_1)}L^{(p_2)}},$$

$$L^{[2]}(p_1,p_2) = 4\mu_{1,L^{(p_1)}L^{(p_2)}} + 2\mu_{2,L^{(p_1)}L^{(p_1)}},$$

$$L^{[3]}(p_1,p_2) = 4\mu_{0,L^{(p_2)}L^{(p_2)}} + 2\mu_{1,L^{(p_1)}L^{(p_2)}}$$

for nonnegative integers $p_1$ and $p_2$, and

$$\hat{\lambda}^2_{p_2|p_1}(\beta) = L^{[1]}(p_1,p_2)\hat{\psi}_\beta(0|0,2) - L^{[2]}(p_1,p_2)\hat{\psi}_\beta(0|2,1) + \mu_{2,L^{(p_2)}L^{(p_2)}}\hat{\psi}_\beta(0|4,0),$$

$$\hat{\kappa}^2_{p_2}(\beta') = 4\mu_{0,L^{(p_2)}L^{(p_2)}}\hat{\psi}_{\beta'}(0|4,0)$$

for bandwidths $\beta$ and $\beta'$. Detailed calculations needed to derive some of equations in the sequel are omitted, but are available from the author.

1. Compute $\widetilde{\mathcal{N}}[6]$ and $\widetilde{\mathcal{D}}[6]$.
2. Compute $\hat{\lambda}^2_{6|7}(\beta_{n1})$ and $\hat{\kappa}^2_6(\beta_{d1})$ for some appropriately chosen bandwidths $\beta_{n1}$ and $\beta_{d1}$, and then compute

$$g_{n1} = \left[\left(\frac{13}{2}\right) \frac{\hat{\lambda}^2_{6|7}(\beta_{n1})}{\mu^2_{2,L}\widetilde{\mathcal{N}}[6]^2}\right]^{1/17} n^{-2/17},$$

$$g_{d1} = \left[\left(\frac{13}{2}\right) \frac{\hat{\kappa}^2_6(\beta_{d1})}{\mu^2_{2,L}\widetilde{\mathcal{D}}[6]^2}\right]^{1/17} n^{-2/17}.$$

3. Compute $\hat{\lambda}^2_{4|5}(\beta_{n2})$ and $\hat{\kappa}^2_4(\beta_{d2})$ for some appropriately chosen bandwidths $\beta_{n2}$ and $\beta_{d2}$, and then compute

$$g_{n2} = \left[\left(\frac{9}{2}\right) \frac{\hat{\lambda}^2_{4|5}(\beta_{n2})}{\mu^2_{2,L}\widehat{\mathcal{N}}_{g_{n1}}[4]^2}\right]^{1/13} n^{-2/13},$$

$$g_{d2} = \left[\left(\frac{9}{2}\right) \frac{\hat{\kappa}^2_4(\beta_{d2})}{\mu^2_{2,L}\widehat{\mathcal{D}}_{g_{d1}}[4]^2}\right]^{1/13} n^{-2/13}.$$

4. Compute $\hat{\lambda}^2_{2|3}(\beta_{n3})$ and $\hat{\kappa}^2_2(\beta_{d3})$ for some appropriately chosen bandwidths $\beta_{n3}$ and $\beta_{d3}$, and then compute

$$g_{n3} = \left[\left(\frac{5}{2}\right)\frac{\hat{\lambda}^2_{2|3}(\beta_{n3})}{\mu^2_{2,L}\widehat{\mathcal{N}}_{g_{n2}}[2]^2}\right]^{1/9} n^{-2/9},$$

$$g_{d3} = \left[\left(\frac{5}{2}\right)\frac{\hat{\kappa}^2_2(\beta_{d3})}{\mu^2_{2,L}\widehat{\mathcal{D}}_{g_{d2}}[2]^2}\right]^{1/9} n^{-2/9}.$$

5. Compute

$$\begin{aligned} g^*_{\mathrm{AMSRE}} = &\left[\frac{5}{2\mu^2_{2,L}\{\widehat{\mathcal{D}}_{g_{d3}}\widehat{\mathcal{N}}_{g_{n2}}[2] - \widehat{\mathcal{N}}_{g_{n3}}\widehat{\mathcal{D}}_{g_{d2}}[2]\}^2}\right]^{1/9} \\ &\times \Big[\widehat{\mathcal{D}}^2_{g_{d3}} L^{[1]}(3,2)\hat{\psi}_{\beta_0}(0|0,2) \\ &\quad - \{\widehat{\mathcal{D}}^2_{g_{d3}} L^{[2]}(3,2) + 2\widehat{\mathcal{N}}_{g_{n3}}\widehat{\mathcal{D}}_{g_{d3}} L^{[3]}(3,2)\}\hat{\psi}_{\beta_0}(0|2,1) \\ &\quad + \{\widehat{\mathcal{D}}^2_{g_{d3}}\mu_{2,L^{(2)}L^{(2)}} + 4\widehat{\mathcal{N}}^2_{g_{n3}}\mu_{0,L^{(2)}L^{(2)}} \\ &\qquad\qquad + 4\widehat{\mathcal{N}}_{g_{n3}}\widehat{\mathcal{D}}_{g_{d3}}\mu_{1,L^{(2)}L^{(3)}}\}\hat{\psi}_{\beta_0}(0|4,0)\Big]^{1/9} \\ &\times n^{-2/9} \end{aligned}$$

for some appropriately chosen bandwidth $\beta_0$.

6. Compute two estimates of $\alpha_o$ defined as

$$\hat{\alpha}^{[2]}_o = \hat{\alpha}_o(g^*_{\mathrm{AMSRE}}) \tag{6.12}$$

and

$$\hat{\alpha}^{[3]}_o = 1 + \frac{1}{2}\frac{\widehat{\mathcal{N}}_{g_{n3}}}{\widehat{\mathcal{D}}_{g_{d3}}}. \tag{6.13}$$

Here $\hat{\alpha}^{[2]}_o$ is based on AMSRE formula (6.11), so that a single bandwidth is included. On the other hand, the two bandwidths included in $\hat{\alpha}^{[3]}_o$ are based on AMSE formulas (6.9) and (6.10), which correspond to the numerator $\mathcal{N}$ and the denominator $\mathcal{D}$, respectively. The bandwidths $\beta_{n1}$, $\beta_{n2}$, $\beta_{n3}$, $\beta_{d1}$, $\beta_{d2}$, $\beta_{d3}$ and $\beta_0$ are all determined using the formula

$$\begin{aligned} &\mathrm{AMSE}\,[a\hat{\psi}_\beta(0|0,2) + b\hat{\psi}_\beta(0|2,1) + c\hat{\psi}_\beta(0|4,0)] \\ &\quad = \frac{\beta^4}{4}\mu^2_{2,L}\{a\psi(2|0,2) + b\psi(2|2,1) + c\psi(2|4,0)\}^2 \\ &\qquad + \frac{2R(L)}{n^2\beta}\int f(x)^2\{aq_2(x)^2 + bq_1(x)^2 q_2(x) + cq_1(x)^4\}^2\,dx \end{aligned}$$

for some constants $a$, $b$ and $c$. This gives the optimal $\beta$ as

$$\beta_{\mathrm{AMSE}} = \left[\frac{2R(L)\,E_f[f(X)\{aq_2(X)^2 + bq_1(X)^2 q_2(X) + cq_1(X)^4\}^2]}{\{a\psi(2|0,2) + b\psi(2|2,1) + c\psi(2|4,0)\}^2}\right]^{1/5} n^{-2/5}.$$

At this stage, estimates of $\psi(0|r,s)$ and $\psi(2|r,s)$ for some pairs $(r,s)$ are needed. These can be provided by kernel estimates of $f$ and $f^{(2)}$ that have bandwidths obtained by the method of Härdle, Marron and Wand (1990). The empirical behavior of $\hat{f}_\alpha$ for $\alpha = \hat{\alpha}_o^{[1]}$, $\hat{\alpha}_o^{[2]}$ and $\hat{\alpha}_o^{[3]}$ is reported in the next section.

**7. Finite sample performance.** Finite sample performance of the proposed density estimators was investigated by Monte Carlo simulation. The first 10 densities (#1–#10) of Marron and Wand (1992), which cover a large variety of realistic density shapes, were used as target densities in this simulation study. In each case 1000 samples of size $n = 500$ were generated. The MISE($h$) value for a given bandwidth $h$ was estimated by the average of these 1000 realizations of (integrated squared error) ISE($h$). To obtain a precise approximation to the minimum MISE, a grid search of the bandwidth was implemented. This was done after an initial screening had provided a suitable $h$ interval that contained the minimum. The Gaussian kernel was used throughout. The estimators compared in this study were $\tilde{f}$ and $\hat{f}_\alpha$ for $\alpha = 0, 1, 2, \alpha_o, \hat{\alpha}_o^{[k]}$, $k = 1, 2, 3$ [see (6.6), (6.12) and (6.13)]. We utilized $g(x, \hat{\theta}) = \phi_{\hat{\sigma}}(x - \hat{\mu})$ for all cases, where $(\hat{\mu}, \hat{\sigma}^2)$ is the MLE of $(\mu, \sigma^2)$. Values of $10^5 \times \min \mathrm{MISE}$ are tabulated in Table 3, where the minimum is taken over $h$. Also tabulated in parentheses for all cases and estimators are $10^5$ times the standard error (SE) of the estimates of MISE($h$) using the bandwidth at which min MISE is obtained.

First we see #1. This case is that $f$ is in the parametric model so that the $O(h^2)$ term of the bias of $\hat{f}_\alpha$ vanishes, as mentioned in Section 3. Therefore, $\alpha_o$ is not defined and the estimation of $\alpha_o$ does not have meaning. Thus, $\hat{f}_\alpha$ for $\alpha = \alpha_o, \hat{\alpha}_o^{[k]}$, $k = 1, 2, 3$, were not simulated for #1 for this reason. For #1 all of $\hat{f}_\alpha$ are significantly better than $\tilde{f}$, and $\hat{f}_1$ is the best.

The tabulated values for $\hat{f}_\alpha$ with $\alpha = \hat{\alpha}_o^{[2]}$ and $\alpha = \hat{\alpha}_o^{[3]}$ in #4 are the median of ISE($h$) for a given $h$ rather than MISE($h$), and the values in parentheses are robust SEs calculated by substituting median absolute deviations. This is because the values of MISE($h$) of these became huge and showed unstable behavior in #4. The instability in #4 can actually be observed, since even the value of robust median ISE($h$) is somewhat large relative to MISE($h$) values of other estimators, and then the value of robust SE in parentheses is also large. We can further observe from #5 that $\hat{f}_\alpha$ with $\alpha = \hat{\alpha}_o^{[3]}$ behaves unstably.

In all cases except #4, #5 and #9, the ideal estimator $\hat{f}_{\alpha_o}$ is the best,

TABLE 3
*The value of estimated* $\min_h \mathrm{MISE}(h)$ $\times 10^5$ *for samples of size* $n = 500$ *from each of the first* 10 *Marron and Wand densities over* 1000 *simulations for* $\tilde{f}$, $\hat{f}_0$ $(= \hat{f}_{\mathrm{HJ}})$, $\hat{f}_1$ $(= \hat{f}_{\mathrm{LL}})$, $\hat{f}_2$ $(= \hat{f}_{\mathrm{HG}})$, $\hat{f}_{\alpha_o}$ *and* $\hat{f}_\alpha$ *with* $\alpha = \hat{\alpha}_o^{[k]}$, $k = 1, 2, 3$. *The standard error* $\times 10^5$ *is given in parentheses for each case*

| | | $\hat{f}_\alpha$ | | | | | | |
|---|---|---|---|---|---|---|---|---|
| **f** | $\tilde{f}$ | $\alpha = 0$ | $\alpha = 1$ | $\alpha = 2$ | $\alpha = \alpha_o$ | $\alpha = \hat{\alpha}_o^{[1]}$ | $\alpha = \hat{\alpha}_o^{[2]}$ | $\alpha = \hat{\alpha}_o^{[3]}$ |
| #1 | 172 | 67 | 62 | 63 | — | — | — | — |
| | (3) | (2) | (1) | (1) | — | — | — | — |
| #2 | 254 | 243 | 196 | 190 | 182 | 227 | 288 | 218 |
| | (4) | (4) | (4) | (4) | (4) | (4) | (9) | (6) |
| #3 | 1,413 | 1,406 | 1,395 | 1,394 | 1,394 | 1,394 | 1,394 | 1,395 |
| | (15) | (15) | (15) | (15) | (15) | (15) | (15) | (15) |
| #4 | 1,372 | 1,296 | 1,290 | 1,286 | 2,734 | 1,288 | 1,440* | 1,523* |
| | (16) | (17) | (17) | (17) | (621) | (17) | (195)† | (213)† |
| #5 | 1,735 | 1,763 | 1,677 | 1,641 | 1,710 | 1,648 | 1,641 | 289,637 |
| | (32) | (32) | (31) | (30) | (28) | (30) | (30) | (5,181) |
| #6 | 244 | 272 | 243 | 234 | 234 | 258 | 234 | 235 |
| | (4) | (4) | (4) | (4) | (4) | (4) | (4) | (4) |
| #7 | 340 | 372 | 340 | 333 | 332 | 336 | 332 | 332 |
| | (5) | (5) | (5) | (5) | (5) | (5) | (5) | (5) |
| #8 | 323 | 361 | 328 | 324 | 321 | 341 | 321 | 324 |
| | (4) | (5) | (5) | (5) | (5) | (5) | (5) | (5) |
| #9 | 296 | 327 | 302 | 297 | 296 | 309 | 296 | 296 |
| | (4) | (4) | (4) | (4) | (4) | (4) | (4) | (4) |
| #10 | 1,126 | 1,139 | 1,125 | 1,124 | 1,123 | 1,135 | 1,124 | 1,124 |
| | (10) | (10) | (10) | (10) | (10) | (10) | (10) | (10) |

Note: The asterisk (∗) designates the minimum of median ISE and the dagger (†) denotes robust SE using median absolute deviation.

which justifies the theory presented in Section 4. In #4 and #5 $\hat{f}_{\alpha_o}$ is not so good because the value $\alpha_o$ is large relative to the other cases as seen in Table 1. It seems that a larger sample is needed for #4 and #5 to confirm the theory presented in Section 4. In addition, good performance of $\hat{f}_{\alpha_o}$ reveals that the estimation of $\alpha_o$ is indeed an important problem. Estimators $\hat{f}_\alpha$ for $\alpha = \hat{\alpha}_o^{[k]}$, $k = 1, 2, 3$, behave well and their differences are small in almost all cases. For $\alpha = \hat{\alpha}_o^{[k]}$, $k = 2, 3$, however, $\hat{f}_\alpha$ were somewhat unstable relative to $\hat{\alpha}_o^{[1]}$ in the sense of the SE, but the bias of these estimators was smaller than that for $\hat{\alpha}_o^{[1]}$.

Some notable insights from Table 3 are as follows. For almost all cases, $\hat{f}_2$ surpasses $\hat{f}_\alpha$ for $\alpha = 0, 1$. Although the degree of improvement is marginal, use of the estimator of $\alpha_o$ yields better performance, which is recognized in #3, #5, #7, #8, #9 and #10. For practical situations, the choices of $\hat{\alpha}_o^{[2]}$ and $\hat{\alpha}_o^{[3]}$ are recommended for densities that are somewhat smooth, but $\hat{\alpha}_o^{[1]}$ is suited for densities that are rather kurtotic.

**8. Supplements.** In this section a number of supplementary results are presented.

8.1. *The integral.* Direct calculation yields

$$\int \hat{f}_\alpha(x)\,dx = 1 + \frac{h^2}{2}\,\mu_{2,K}\,\frac{1}{n}\sum_{i=1}^{n}\left[\frac{\{g(X_i,\hat{\theta})^{\alpha-1}\}''}{g(X_i,\hat{\theta})^{\alpha-1}} - \frac{\{g(X_i,\hat{\theta})^{2-\alpha}\}''}{g(X_i,\hat{\theta})^{2-\alpha}}\right] + O(h^4)$$

as $h \to 0$. In particular, when we adopt the Gaussian density $g(x,\hat{\theta}) = \phi_{\hat{\sigma}}(x-\hat{\mu}) = \phi_s(x-\bar{X})$ as an initial parametric start, where $\bar{X}$ and $s^2$ are, respectively, the sample mean and the sample variance, we have

$$\begin{aligned}\int \hat{f}_\alpha(x)\,dx &= 1 + \frac{h^2}{2}\mu_{2,K}\left(\frac{2\alpha-3}{s^2}\right)\frac{1}{n}\sum_{i=1}^{n}\left\{\left(\frac{X_i-\bar{X}}{s}\right)^2 - 1\right\} + O(h^4)\\ &= 1 + O(h^4)\end{aligned}$$

as $h \to 0$.

8.2. *Computational remark.* The practical expression for $\hat{f}_\alpha$ depends on the choices of the kernel $K$ and the initial parametric model $g(x,\theta)$. Thus, the general features required for practical calculation are not pursued here. However, derivation of the expression for the case in which the Gaussian kernel and model are adopted appears to be useful. Now define that

$$\mathcal{I}(\alpha) = \int K_h(t-x)g(t,\hat{\theta})^{-\alpha}\,dt$$

for $K(t) = \phi(t)$ and $g(t,\hat{\theta}) = \phi_{\hat{\sigma}}(t-\hat{\mu})$. Direct calculations give

$$\mathcal{I}(\alpha) = \frac{(\sqrt{2\pi}\,)^\alpha \hat{\sigma}^{\alpha+1}}{\sqrt{\hat{\sigma}^2-\alpha h^2}} \exp\left[\frac{\alpha(x-\hat{\mu})^2}{2(\hat{\sigma}^2-\alpha h^2)}\right],$$

provided that $\hat{\sigma}^2 - \alpha h^2 > 0$. Using this notation, we have, for the case of Gaussian kernel and model,

$$\begin{aligned}\hat{f}_\alpha(x) &= \frac{(\sqrt{2\pi}\,)^{\alpha-3}\hat{\sigma}^{\alpha-2}}{nh\mathcal{I}(\alpha-2)}\\ &\quad\times\sum_{i=1}^{n}\exp\left\{-\frac{(x-\hat{\mu})^2}{2\hat{\sigma}^2} - \frac{(X_i-x)^2}{2h^2} - (1-\alpha)\frac{(X_i-\hat{\mu})^2}{2\hat{\sigma}^2}\right\}\end{aligned}$$

for $\hat{\sigma}^2 - (\alpha-2)h^2 > 0$.

8.3. *Choosing the bandwidth.* From Section 4 we see that the bandwidth $h$ that minimizes the AMISE for $\hat{f}_\alpha$ is

$$h(\alpha) = \left\{ \frac{R(K)}{\mu_{2,K}^2} \right\}^{1/5} \mathcal{R}(\hat{f}_\alpha)^{-1/5} n^{-1/5}$$

for a fixed $\alpha$, and the resultant minimum value of the AMISE is

$$\tfrac{5}{4}\{\mu_{2,K} R(K)^2\}^{2/5} \mathcal{R}(\hat{f}_\alpha)^{1/5} n^{-4/5}.$$

Proposition 1 reveals that we can further reduce this by using $\alpha = \alpha_o$ in (4.10). Thus, the best choice for the bandwidth $h$ is

$$h_o = h(\alpha_o) = \left\{ \frac{R(K)}{\mu_{2,K}^2} \right\}^{1/5} \mathcal{R}(\hat{f}_{\alpha_o})^{-1/5} n^{-1/5}.$$

Here, we propose a method to choose $h$ which is a variant of that discussed in Hjort and Glad (1995). Recall the analogy presented in Section 6.1, and consider $\bar{h}$ in (6.5) and $\hat{\alpha}_o^{[1]}$ in (6.6). Further, we consider a bias-adjusted version of $\mathcal{R}(\hat{f}_\alpha)$ given as

$$\mathcal{R}^\dagger(\alpha, h) = \frac{n}{n-1} \left\{ \hat{\mathcal{R}}(\alpha, h) - \frac{R(K'')}{nh^5} \right\},$$

where

$$\hat{\mathcal{R}}(\alpha, h) = \hat{c}_1(h)\alpha^2 - 2\hat{c}_2(h)\alpha + \hat{c}_3(h)$$

and

$$\hat{c}_3(h) = \int \{\hat{b}_1(x; h) + \hat{b}_2(x; h)\}^2 \, dx.$$

Here $\bar{h}$ in (6.5) is seen as an initial bandwidth. Then we calculate the final bandwidth as

$$\hat{h} = \left\{ \frac{R(K)}{\mu_{2,K}^2} \right\}^{1/5} \mathcal{R}^\dagger(\hat{\alpha}_o^{[1]}, \bar{h})^{-1/5} n^{-1/5}.$$

The theoretical performance of this $\hat{h}$ is not pursued here. However, we have an empirical suggestion based on application to some artificial data that $\hat{h}$ is not as stable as $\bar{h}$.

**9. Proofs.** In this section the proofs of theoretical results are presented. First, we prepare the following lemma which can be proved by Taylor expansion.

LEMMA 1. *Let* $g_0(x) = g(x, \theta_0)$ *and*

$$f_\alpha^*(x) = g_0(x) \frac{n^{-1} \sum_{i=1}^n K_h(X_i - x) g_0(X_i)^{1-\alpha}}{\int K_h(t - x) g_0(t)^{2-\alpha} \, dt}. \tag{9.1}$$

*Then as* $n \to \infty$, $h \to 0$,

$$\operatorname{Bias} f^*_\alpha(x) = \frac{h^2}{2}\mu_{2,K}\left[\frac{(g_0(x)^{1-\alpha}f(x))''}{g_0(x)^{1-\alpha}} - \frac{f(x)(g_0(x)^{2-\alpha})''}{g_0(x)^{2-\alpha}}\right] + O(h^4),$$

$$\operatorname{Var} f^*_\alpha(x) = \frac{R(K)}{nh}f(x) - \frac{f(x)^2}{n} + O\left(\frac{h}{n}\right).$$

Proof of Proposition 1. The result is straightforwardly obtained from the quadratic expression of $\mathcal{R}(\hat{f}_\alpha)$ in (4.9). □

Proof of Theorem 1. Define

$$u_0(x) = \frac{\partial}{\partial\theta}\log g(x,\theta_0),$$

$$U_0(x) = \frac{\partial^2}{\partial\theta\,\partial\theta^T}\log g(x,\theta_0).$$

Using Taylor expansions, we can expand $\hat{f}_\alpha$ as

$$\hat{f}_\alpha(x) = f^*_\alpha(x) + (\hat\theta - \theta_0)^T\bar{B}_n(x) + \tfrac{1}{2}(\hat\theta-\theta_0)^T\bar{C}_n(x)(\hat\theta-\theta_0) + o_p(n^{-1}),$$

where $f^*_\alpha$ is given as in (9.1),

$$\bar{B}_n(x) = \frac{1}{n}\sum_{i=1}^n B_i(x),$$

$$\bar{C}_n(x) = \frac{1}{n}\sum_{i=1}^n C_i(x),$$

$$B_i(x) = K_h(X_i - x)g_0(X_i)^{1-\alpha}\frac{g_0(x)}{\eta_0(x)} \times\left[(1-\alpha)u_0(X_i) - \frac{(2-\alpha)}{\eta_0(x)}\eta_1(x) + u_0(x)\right],$$

$$\begin{aligned} C_i(x) = {} & K_h(X_i - x)g_0(X_i)^{1-\alpha}\frac{g_0(x)}{\eta_0(x)} \\ & \times\left[-\frac{2(1-\alpha)(2-\alpha)}{\eta_0(x)}\eta_1(x)u_0(X_i)^T + 2(1-\alpha)u_0(x)u_0(X_i)^T\right. \\ & + (1-\alpha)\{U_0(X_i) + (1-\alpha)u_0(X_i)u_0(X_i)^T\} \\ & - \frac{2(2-\alpha)}{\eta_0(x)}u_0(x)\eta_1(x)^T + \{U_0(x) + u_0(x)u_0(x)^T\} \\ & \left. + \frac{2(2-\alpha)}{\eta_0(x)^2}\left\{(2-\alpha)\eta_1(x)\eta_1(x)^T - \frac{1}{2}\eta_0(x)\eta_2(x)\right\}\right], \end{aligned}$$

where

$$\eta_0(x) = \int K_h(t-x) g_0(t)^{2-\alpha}\, dt,$$

$$\eta_1(x) = \int K_h(t-x) u_0(t) g_0(t)^{2-\alpha}\, dt,$$

$$\eta_2(x) = \int K_h(t-x)\{U_0(t) + (2-\alpha)u_0(t)u_0(t)^T\} g_0(t)^{2-\alpha}\, dt.$$

Through (3.1) and the average representations above, we have

$$E[(\hat{\theta}-\theta_0)^T \bar{B}_n(x)] = O\left(\frac{h^2}{n} + \frac{1}{n^2}\right),$$

$$E[(\hat{\theta}-\theta_0)^T \bar{C}_n(x)(\hat{\theta}-\theta_0)] = O\left(\frac{h^2}{n} + \frac{1}{n^2}\right),$$

using the fact that $I_i = I(X_i)$ has mean zero. Since the bias term of $f^*_\alpha$ in (9.1) was already given in Lemma 1, the bias expression of $\hat{f}_\alpha$ is confirmed.

Next we consider variance. The variance of $f^*_\alpha$ was obtained in Lemma 1. By using the average representation (3.1), we have, after somewhat lengthy calculations, that

$$\operatorname{Var}[(\hat{\theta}-\theta_0)^T \bar{B}_n(x)] = O\left(\frac{h^4}{n} + \frac{1}{n^2}\right),$$

$$\operatorname{Var}[(\hat{\theta}-\theta_0)^T \bar{C}_n(x)(\hat{\theta}-\theta_0)] = O\left(\frac{h^4}{n^2}\right),$$

$$\operatorname{Cov}[f^*_\alpha(x), (\hat{\theta}-\theta_0)^T \bar{B}_n(x)] = O\left(\frac{h^2}{n} + \frac{1}{n^2}\right),$$

from which the necessary variance expression is derived. □

Proof of Theorem 2. Direct calculation yields that

$$\begin{aligned}
&\operatorname{MSRE}[\hat{\alpha}_o(g)] \\
&\quad = E\left[\left\{\frac{\hat{\alpha}_o(g)}{\alpha_o} - 1\right\}^2\right] \\
&\quad = \frac{\mathcal{D}^2}{4\mathcal{N}^2} E\left[\left\{\frac{\mathcal{D}(\widehat{\mathcal{N}}_g - \mathcal{N}) - \mathcal{N}(\widehat{\mathcal{D}}_g - \mathcal{D})}{\mathcal{D}^2 + \mathcal{D}(\widehat{\mathcal{D}}_g - \mathcal{D})}\right\}^2\right] \\
&\quad = \frac{1}{4}\left[\frac{\operatorname{MSE}[\widehat{\mathcal{N}}_g]}{\mathcal{N}^2} + \frac{\operatorname{MSE}[\widehat{\mathcal{D}}_g]}{\mathcal{D}^2} - 2\frac{E\{(\widehat{\mathcal{D}}_g - \mathcal{D})(\widehat{\mathcal{N}}_g - \mathcal{N})\}}{\mathcal{N}\mathcal{D}}\right] \\
&\qquad + O_{n,g},
\end{aligned}$$

where $O_{n,g}$ is a negligible higher-order term. Hence it suffices to show (6.9) and (6.10), and to evaluate the cross term $E\{(\widehat{\mathcal{D}}_g - \mathcal{D})(\widehat{\mathcal{N}}_g - \mathcal{N})\}$ for checking (6.11). However, only the proof of (6.9) is presented here since the other equations can be obtained in the same manner. We therefore focus on $\widehat{\mathcal{N}}_g$. Then it follows that

$$
\begin{aligned}
&\mathrm{MSE}[\widehat{\mathcal{N}}_g] \\
&\quad = \mathrm{MSE}[\hat{\psi}_g(0|2,1)] + \mathrm{MSE}[\hat{\psi}_g(3|1,0)] \\
&\qquad + \mathrm{MSE}[\hat{\psi}_g(2|0,1)] + \mathrm{MSE}[\hat{\psi}_g(1|1,1)] \\
&\qquad - 2E[\hat{\mu}_g(0|2,1)\hat{\mu}_g(3|1,0)] - 2E[\hat{\mu}_g(0|2,1)\hat{\mu}_g(2|0,1)] \\
&\qquad - 2E[\hat{\mu}_g(0|2,1)\hat{\mu}_g(1|1,1)] + 2E[\hat{\mu}_g(3|1,0)\hat{\mu}_g(2|0,1)] \\
&\qquad + 2E[\hat{\mu}_g(3|1,0)\hat{\mu}_g(1|1,1)] + 2E[\hat{\mu}_g(2|0,1)\hat{\mu}_g(1|1,1)],
\end{aligned} \tag{9.2}
$$

where $\hat{\mu}_g(p|r,s) = \hat{\psi}_g(p|r,s) - \psi(p|r,s)$. Therefore, the proof is further reduced to evaluation of $\mathrm{MSE}[\hat{\psi}_g(p|r,s)]$ and $E[\hat{\mu}_g(p_1|r_1,s_1)\hat{\mu}_g(p_2|r_2,s_2)]$ for nonnegative integer triplets $(p|r,s)$, $(p_1|r_1,s_1)$ and $(p_2|r_2,s_2)$. To accomplish the proof, the following four lemmas are needed. The proofs of all four lemmas are omitted. Details are available from the author.

Let us define

$$
\psi_g^*(p|r,s) = \frac{1}{n(n-1)} \sum_{i\neq j} q_1(X_i)^r q_2(X_i)^s L_g^{(p)}(X_i - X_j). \tag{9.3}
$$

Performance of $\hat{\psi}_g(p|r,s)$ is dominated by the performance of $\psi_g^*(p|r,s)$. The following lemma is concerned with $\psi_g^*(p|r,s)$.

LEMMA 2. *Let $\psi_g^*(p|r,s)$ be as given in (9.3). Then, as $n\to\infty$, $g\to 0$,*

$$
\begin{aligned}
&\mathrm{MSE}[\psi_g^*(p|r,s)] \\
&\quad = \mathrm{Bias}[\psi_g^*(p|r,s)]^2 + \mathrm{Var}[\psi_g^*(p|r,s)] \\
&\quad = \frac{g^4}{4}\mu_{2,L}^2 \psi(p+2|2r,s)^2 + \frac{2R(L^{(p)})}{n^2 g^{2p+1}} \psi(0|2r,2s) \\
&\qquad + \frac{1}{n}\left[\int f(x)\{w(x) f^{(p)}(x) + \{w\cdot f\}^{(p)}(x)\}^2\,dx - 4E[\psi_g^*(p|r,s)]^2\right] \\
&\qquad + o(n^{-1} + n^{-2}g^{-2p-1}),
\end{aligned}
$$

*for $p$ even and*

$$
\begin{aligned}
&\mathrm{MSE}[\psi_g^*(p|r,s)] \\
&\quad = \frac{g^4}{4}\mu_{2,L}^2 \psi(p+2|r,s)^2 \\
&\qquad \times \frac{\mu_{2,\{L^{(p)}\}^2}}{2n^2 g^{2p-1}} \int f(x)\{\{w^2\cdot f\}^{(2)}(x) - w(x)\{w\cdot f\}^{(2)}(x)\}\,dx
\end{aligned}
$$

$$+\frac{1}{n}\Big[\int f(x)\{w(x)f^{(p)}(x)-\{w\cdot f\}^{(p)}(x)\}^2\,dx-4E[\psi_g^*(p|r,s)]^2\Big]$$
$$+o(n^{-1}+n^{-2}g^{-2p+1})$$

*for $p$ odd, where*

$$w(x)=q_1(x)^r q_2(x)^s.$$

The notation

$$w_{r,s}(x)=q_1(x)^r q_2(x)^s,\qquad \phi_p(x)=L_g^{(p)}(x)$$

is used in the next lemma.

LEMMA 3. *As $n\to\infty$, $g\to 0$, we have*

$$\begin{aligned}
&E[\psi_g^*(p_1|r_1,s_1)\psi_g^*(p_2|r_2,s_2)]\\
&\quad=E[w_{r_1,s_1}(X_1)\phi_{p_1}(X_1-X_2)]E[w_{r_2,s_2}(X_1)\phi_{p_2}(X_1-X_2)]\\
&\qquad+\frac{\mu_{1,L^{(p_1)}L^{(p_2)}}}{n^2g^{p_1+p_2}}\\
&\qquad\quad\times\int f[\{w_{r_1+r_2,s_1+s_2}\cdot f\}^{(1)}+(-1)^{p_2}w_{r_2,s_2}\{w_{r_1,s_1}\cdot f\}^{(1)}](x)\,dx\\
&\qquad+\frac{1}{n}\Big[\int f\{w_{r_1,s_1}f^{(p_1)}+(-1)^{p_1}\{w_{r_1,s_1}\cdot f\}^{(p_1)}\}\\
&\qquad\qquad\times\{w_{r_2,s_2}f^{(p_2)}+(-1)^{p_2}\{w_{r_2,s_2}\cdot f\}^{(p_2)}\}(x)\,dx\\
&\qquad\quad-4E[w_{r_1,s_1}(X_1)\phi_{p_1}(X_1-X_2)]E[w_{r_2,s_2}(X_1)\phi_{p_2}(X_1-X_2)]\Big]\\
&\qquad+o(n^{-1}+n^{-2}g^{-p_1-p_2})
\end{aligned}$$

*for $p_1+p_2$ odd, and*

$$\begin{aligned}
&E[\psi_g^*(p_1|r_1,s_1)\psi_g^*(p_2|r_2,s_2)]\\
&\quad=E[w_{r_1,s_1}(X_1)\phi_{p_1}(X_1-X_2)]E[w_{r_2,s_2}(X_1)\phi_{p_2}(X_1-X_2)]\\
&\qquad+\frac{2}{n^2g^{p_1+p_2+1}}\psi(0|r_1+r_2,s_1+s_2)\int L^{(p_1)}(z)L^{(p_2)}(z)\,dz\\
&\qquad+\frac{1}{n}\Big[\int f\{w_{r_1,s_1}f^{(p_1)}+(-1)^{p_1}\{w_{r_1,s_1}\cdot f\}^{(p_1)}\}\\
&\qquad\qquad\times\{w_{r_2,s_2}f^{(p_2)}+(-1)^{p_2}\{w_{r_2,s_2}\cdot f\}^{(p_2)}\}(x)\,dx\\
&\qquad\quad-4E[w_{r_1,s_1}(X_1)\phi_{p_1}(X_1-X_2)]E[w_{r_2,s_2}(X_1)\phi_{p_2}(X_1-X_2)]\Big]\\
&\qquad+o(n^{-1}+n^{-2}g^{-p_1-p_2-1})
\end{aligned}$$

*for $p_1+p_2$ even, with both $p_1$ and $p_2$ being even, and*

$$\begin{aligned}
&E[\psi_g^*(p_1|r_1,s_1)\psi_g^*(p_2|r_2,s_2)]\\
&\quad = E[w_{r_1,s_1}(X_1)\phi_{p_1}(X_1-X_2)]E[w_{r_2,s_2}(X_1)\phi_{p_2}(X_1-X_2)]\\
&\qquad + \frac{\mu_{2,L^{(p_1)}L^{(p_2)}}}{2n^2g^{p_1+p_2-1}}\int f[\{w_{r_1+r_2,s_1,s_2}\cdot f\}^{(2)} - w_{r_2,s_2}\{w_{r_1,s_1}\cdot f\}^{(2)}](x)\,dx\\
&\qquad + \frac{1}{n}\biggl[\int f\{w_{r_1,s_1}f^{(p_1)} + (-1)^{p_1}\{w_{r_1,s_1}\cdot f\}^{(p_1)}\}\\
&\qquad\qquad \times\{w_{r_2,s_2}f^{(p_2)} + (-1)^{p_2}\{w_{r_2,s_2}\cdot f\}^{(p_2)}\}(x)\,dx\\
&\qquad\quad - 4E[w_{r_1,s_1}(X_1)\phi_{p_1}(X_1-X_2)]E[w_{r_2,s_2}(X_1)\phi_{p_2}(X_1-X_2)]\biggr]\\
&\qquad + o(n^{-1}+n^{-2}g^{-p_1-p_2-1})
\end{aligned}$$

*for both $p_1$ and $p_2$ being odd.*

Hereafter, we adopt the notation

$$\begin{aligned}
A_n &= \frac{1}{n(n-1)}\sum_{i\neq j} L_g^{(p)}(X_i-X_j)v(X_i),\\
B_n &= \frac{1}{n(n-1)}\sum_{i\neq j} L_g^{(p)}(X_i-X_j)W(X_i),\\
v(x) &= \frac{\partial}{\partial\theta}\left\{\frac{g'(x,\theta)}{g(x,\theta)}\right\}^r\left\{\frac{g''(x,\theta)}{g(x,\theta)}\right\}^s\bigg|_{\theta=\theta_0},\\
W(x) &= \frac{\partial^2}{\partial\theta\,\partial\theta^T}\left\{\frac{g'(x,\theta)}{g(x,\theta)}\right\}^r\left\{\frac{g''(x,\theta)}{g(x,\theta)}\right\}^s\bigg|_{\theta=\theta_0}.
\end{aligned}$$

The behavior of $\hat\psi_g(p|r,s)$ is summarized in the next lemma.

LEMMA 4. *As $n\to\infty$, $g\to 0$, we have*

$$\begin{aligned}
&\mathrm{MSE}[\hat\psi_g(p|r,s)]\\
&\quad = \frac{g^4}{4}\mu_{2,L}^2\psi(p+2|r,s)^2\\
&\qquad + \frac{\mu_{2,\{L^{(p)}\}^2}}{2n^2g^{2p-1}}\int f(x)\{\{w^2\cdot f\}^{(2)}(x) + (-1)^p w(x)\{w\cdot f\}^{(2)}(x)\}\,dx\\
&\qquad + \frac{1}{n}\biggl[\int f(x)\{w(x)f^{(p)}(x) + (-1)^p\{w\cdot f\}^{(p)}(x)\}^2\,dx\\
&\qquad\quad - 4E[\psi_g^*(p|r,s)]^2 + E[A_n]^T\Sigma_I E[A_n]
\end{aligned}$$

$$+ 2\{E[w(X_1)\phi(X_1 - X_2)(I_1 + I_2)]\}^T E[A_n]\Big]$$

$$+ o(n^{-1} + n^{-2}g^{-2p+1})$$

*for p odd and*

$$\begin{aligned}
&\mathrm{MSE}[\hat{\psi}_g(p|r,s)] \\
&= \frac{g^4}{4}\mu_{2,L}^2 \psi(p+2|r,s)^2 + \frac{2}{n^2 g^{2p+1}} R(L^{(p)})\psi(0|2r,2s) \\
&\quad + \frac{1}{n}\Big[\int f(x)\{w(x)f^{(p)}(x) + (-1)^p\{w \cdot f\}^{(p)}(x)\}^2\,dx \\
&\quad\quad - 4E[\psi_g^*(p|r,s)]^2 + E[A_n]^T \Sigma_I E[A_n] \\
&\quad\quad\quad + 2\{E[w(X_1)\phi(X_1 - X_2)(I_1 + I_2)]\}^T E[A_n]\Big] \\
&\quad + o(n^{-1} + n^{-2}g^{-2p-1})
\end{aligned}$$

*for p even.*

Lemma 5. *As $n \to \infty$, $g \to 0$, we have*

$$\begin{aligned}
&E[\hat{\mu}_g(p_1|r_1,s_1)\hat{\mu}_g(p_2|r_2,s_2)] \\
&= \frac{g^4}{4}\mu_{2,L}^2 \psi(p_1+2|r_1,s_1)\psi(p_2+2|r_2,s_2) \\
&\quad + \frac{\mu_{1,L^{(p_1)}L^{(p_2)}}}{n^2 g^{p_1+p_2}} \\
&\quad\quad \times \int f[\{w_{r_1+r_2,s_1+s_2} \cdot f\}^{(1)} + (-1)^{p_2} w_{r_2,s_2}\{w_{r_1,s_1} \cdot f\}^{(1)}](x)\,dx \\
&\quad + \frac{1}{n}\Big[\int f[w_{r_1,s_1} \cdot f^{(p_1)} + (-1)^{p_1}\{w_{r_1,s_1} \cdot f\}^{(p_1)}] \\
&\quad\quad\quad \times [w_{r_2,s_2} \cdot f^{(p_2)} + (-1)^{p_2}\{w_{r_2,s_2} \cdot f\}^{(p_2)}](x)\,dx \\
&\quad\quad - 4E[w_{r_1,s_1}(X_1)\phi_{p_1}(X_1 - X_2)]E[w_{r_2,s_2}(X_1)\phi_{p_2}(X_1 - X_2)] \\
&\quad\quad + E[w_{r_1,s_1}(X_1)\phi_{p_1}(X_1 - X_2)(I_1 + I_2)]^T E[A_n(p_2|r_2,s_2)] \\
&\quad\quad + E[w_{r_2,s_2}(X_1)\phi_{p_2}(X_1 - X_2)(I_1 + I_2)]^T E[A_n(p_1|r_1,s_1)] \\
&\quad\quad\quad + E[A_n(p_1|r_1,s_1)]^T \Sigma_I E[A_n(p_2|r_2,s_2)]\Big] \\
&\quad + o(n^{-1} + n^{-2}g^{-p_1-p_2})
\end{aligned}$$

*for $p_1 + p_2$ odd and*

$$
\begin{aligned}
&E[\hat{\mu}_g(p_1|r_1, s_1)\hat{\mu}_g(p_2|r_2, s_2)] \\
&\quad = \frac{g^4}{4}\mu_{2,L}^2 \psi(p_1 + 2|r_1, s_1)\psi(p_2 + 2|r_2, s_2) \\
&\qquad + \frac{2\mu_{0,L^{(p_1)}L^{(p_2)}}}{n^2 g^{p_1+p_2+1}}\psi(0|r_1 + r_2, s_1 + s_2) \\
&\qquad + \frac{1}{n}\Big[\int f[w_{r_1,s_1}\cdot f^{(p_1)} + (-1)^{p_1}\{w_{r_1,s_1}\cdot f\}^{(p_1)}] \\
&\qquad\qquad \times [w_{r_2,s_2}\cdot f^{(p_2)} + (-1)^{p_2}\{w_{r_2,s_2}\cdot f\}^{(p_2)}](x)\,dx \\
&\qquad\quad - 4E[w_{r_1,s_1}(X_1)\phi_{p_1}(X_1 - X_2)]E[w_{r_2,s_2}(X_1)\phi_{p_2}(X_1 - X_2)] \\
&\qquad\quad + E[w_{r_1,s_1}(X_1)\phi_{p_1}(X_1 - X_2)(I_1 + I_2)]^T E[A_n(p_2|r_2, s_2)] \\
&\qquad\quad + E[w_{r_2,s_2}(X_1)\phi_{p_2}(X_1 - X_2)(I_1 + I_2)]^T E[A_n(p_1|r_1, s_1)] \\
&\qquad\qquad\qquad + E[A_n(p_1|r_1, s_1)]^T \Sigma_I E[A_n(p_2|r_2, s_2)]\Big] \\
&\qquad + o(n^{-1} + n^{-2}g^{-p_1-p_2-1})
\end{aligned}
$$

*for $p_1$ even and $p_2$ even, and*

$$
\begin{aligned}
&E[\hat{\mu}_g(p_1|r_1, s_1)\hat{\mu}_g(p_2|r_2, s_2)] \\
&\quad = \frac{g^4}{4}\mu_{2,L}^2 \psi(p_1 + 2|r_1, s_1)\psi(p_2 + 2|r_2, s_2) \\
&\qquad + \frac{\mu_{2,L^{(p_1)}L^{(p_2)}}}{2n^2 g^{p_1+p_2-1}}\int f[\{w_{r_1+r_2,s_1+s_2}\cdot f\}^{(2)} - w_{r_2,s_2}\{w_{r_1,s_1}\cdot f\}^{(2)}](x)\,dx \\
&\qquad + \frac{1}{n}\Big[\int f[w_{r_1,s_1}\cdot f^{(p_1)} + (-1)^{p_1}\{w_{r_1,s_1}\cdot f\}^{(p_1)}] \\
&\qquad\qquad \times [w_{r_2,s_2}\cdot f^{(p_2)} + (-1)^{p_2}\{w_{r_2,s_2}\cdot f\}^{(p_2)}](x)\,dx \\
&\qquad\quad - 4E[w_{r_1,s_1}(X_1)\phi_{p_1}(X_1 - X_2)]E[w_{r_2,s_2}(X_1)\phi_{p_2}(X_1 - X_2)] \\
&\qquad\quad + E[w_{r_1,s_1}(X_1)\phi_{p_1}(X_1 - X_2)(I_1 + I_2)]^T E[A_n(p_2|r_2, s_2)] \\
&\qquad\quad + E[w_{r_2,s_2}(X_1)\phi_{p_2}(X_1 - X_2)(I_1 + I_2)]^T E[A_n(p_1|r_1, s_1)] \\
&\qquad\qquad\qquad + E[A_n(p_1|r_1, s_1)]^T \Sigma_I E[A_n(p_2|r_2, s_2)]\Big] \\
&\qquad + o(n^{-1} + n^{-2}g^{-p_1-p_2+1})
\end{aligned}
$$

*for $p_1$ odd and $p_2$ odd.*

Proof of Theorem 2 (continued). By applying Lemmas 4 and 5 to (9.2) and rearranging, the MSE expression of $\widehat{\mathcal{N}}_g$ is obtained. This completes the proof. □

**Acknowledgments.** I thank an Associate Editor and a referee for useful comments and suggestions which yielded substantial improvements in the article.

Department of Mathematics
Faculty of Science and Engineering
Shimane University
1060 Nishi-Kawatsu
Matsue, 690-8504
Japan
e-mail: naito@math.shimane-u.ac.jp